\documentstyle[12pt,twoside, epsf]{article}

\newenvironment{pf}{{\bf Proof.}}{\hfill $\Box$ \medskip}

\newtheorem{th}{Theorem}[section]
\newtheorem{lem}[th]{Lemma}
\newtheorem{cor}[th]{Corollary}

\newtheorem{defn}[th]{Definition}

\newtheorem{rem}[th]{Remark}        
\newtheorem{rems}[th]{Remarks}        
\newtheorem{note}[th]{Note}

\font\sc=cmcsc10

\mathchardef\twid="1218

\def\d{\hbox{\spec \char'017\kern 0.05em}} 
\def\text{\hbox}

\let\\=\cr

\def\Chi{\hbox{\raise0.5ex\hbox{$\chi$}}}
\def\C{{\cal C}}
\def\B{{\cal B}}

\def\picill#1by#2(#3)
{\vbox to #2
{\hrule width #1 height 0pt depth 0pt
\vfill\epsffile{#3}}}

\let \ttorg \tt \def \tt{\ttorg \obeyspaces}

\begin{document}
\pagestyle{myheadings}
\markboth{Lambropoulou--Rourke}{Markov's theorem in 3--manifolds}

\title{Markov's theorem in 3--manifolds}
\author{\sc Sofia Lambropoulou and Colin Rourke}
\date{}

\maketitle

\begin{abstract}
In this paper we first give a one-move version of Markov's braid theorem for knot isotopy in
 $S^3$ that sharpens the classical theorem. Then  a relative version of Markov's theorem
concerning a fixed braided portion in the knot. We also prove an analogue of Markov's theorem for
knot isotopy in knot complements. Finally we extend this last result to prove a Markov theorem for
links in an arbitrary orientable 3--manifold. \end{abstract}

\section{Overview}

  According to Birman \cite{B}, Markov \cite{Ma} originally stated his braid equivalence
theorem using {\it three} braid moves; later  there was another brief announcement of
an improved version of Markov's theorem by Weinberg \cite{W} consisting of the {\it two} well-known
braid-equivalence moves: conjugation in the braid groups and the `stabilizing' or `Markov' moves
($M$--moves). Our first main result is a {\it one-move} Markov theorem which we now state. 

An {\it $L$--move} on a braid consists of cutting one arc of
the braid open and splicing into the broken strand new strands to the top
and bottom, both either {\it under} or {\it over} the rest of the
braid:

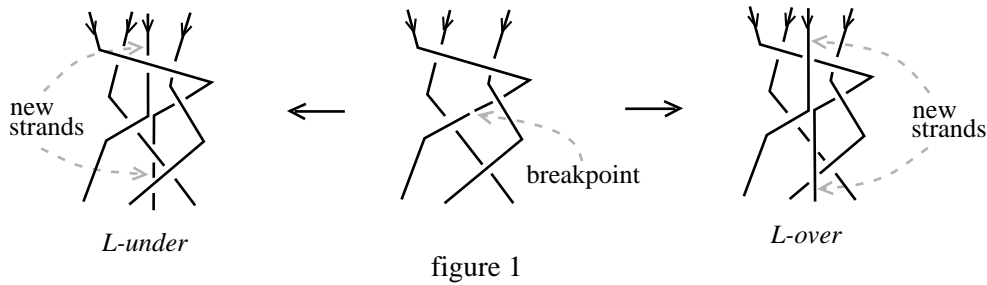
\bigbreak

$$\vbox{\picill5inby1.4in(fig1)  }$$

\noindent $L$--moves and isotopy generate an equivalence relation on
braids called {\it$L$--equivalence}. In \S 4 we prove that
$L$--equivalence classes of braids are in bijective correspondence
with isotopy classes of oriented links in $S^3$, where the bijection is
induced by `closing' the braid to form a link. As a consequence,
$L$--equivalence is the same as the usual Markov equivalence and thus the classical Markov 
theorem is sharpened. The proofs are based on a canonical
process for turning a combinatorial oriented link diagram in the plane (with a
little extra structure) into an open braid. Our braiding as well as the $L$--moves are  based on
the building blocks of combinatorial isotopy, the {\it triangle moves} or {\it $\Delta$--moves}.
This makes the proof
 conceptually very simple\footnote{Our braiding was first used in \cite{La0} and then in 
\cite{La} and in \cite{LR} (preliminary version of our results). The first complete proof of the
classical theorem has  been given by J.~Birman \cite{B}, and other proofs by H.~Morton \cite{M},
D.~Bennequin \cite{Ben} and P. Traczyk \cite{Tr}.}. Our braiding operation is essentially the
same as the operation using a `saw-tooth' given by J. Birman in \cite{B}. We use the point at
infinity as the reference point for braiding and so a saw-tooth becomes a pair of vertical lines
that meet at infinity. The change of reference point to infinity makes the proof of Markov
theorem easier because there are very few ways that saw-teeth can obstruct each other.

Moreover, the local nature of our proof
allows us to formulate and prove first the relative version of Markov theorem, another new
result stating that, if two isotopic links contain the same braided portion then any two
corresponding braids which both contain that braided portion differ by a sequence of $L$--moves that
do not affect the braided portion.  We then prove a second relative version of the theorem for
links which contain a common closed braid and deduce  an analogue of Markov  theorem for knot
isotopy in complements of knots/links, which we finally extend to an analogue of  Markov theorem
for  links  in  arbitrary closed 3--manifolds. 
\bigbreak
 More precisely, let  $M$ be a closed connected orientable 3--manifold.  We can assume that $M$ is
given by surgery on the  closure of a braid $B$ in $S^3$, thus $M$ is represented by $B$.
Furthermore, any oriented link in $M$ can be represented as the closure of a further braid $\beta$
such that $B\cup\beta$  is also a braid and we shall call it  {\it mixed braid}. Then the mixed
braid equivalence in $S^3$ that reflects isotopy in $M$ is generated firstly by $L$--moves, as
above, but performed only on strings of $\beta$ and secondly by {\it braid band moves} or {\it
b.b.--moves}. These are moves that  reflect the sliding of a part of a link across the 2--disc
bounded by the specified longitude of a surgery component of the closure of $B$.  I.e. suppose a
string of $\beta$ is adjacent to one of  $B$, then a slide of the first over the second (with a
half-twisted band) replaces $\beta$ by a braid with one or more extra strings parallel to the
strings of $B$ which form the appropriate component in the closure:

\bigbreak
$$\vbox{\picill4inby1.4in(fig2)  }$$

Most  results are based on material in \cite{La} worked at Warwick University under the advice of
the second author. However, they are sharpened versions and some parts of the proofs are improved
considerably.  
\bigbreak
\noindent {\bf Acknowledgements} \  The second author would like to acknowledge a conversation
with Dale Rolfsen which was crucial in formulating the one-move Markov theorem. The first author
wishes  to acknowledge discussions with Jozef
Przytycki, Cameron Gordon and David Yetter on the set--up and some technical details, in Warsaw
'95, and also to thank the E.U. for financial support and the SFB 170 in G\"{o}ttingen for
hospitality and financial support.

\section{Diagrams, isotopy and the $L$-moves}

A knot is a special case of a link, so from  now on we shall be referring to both knots and links as
`links'. Throughout the first part of this paper we shall be working in a combinatorial
setting. In particular, we consider oriented {\it links/geometric braids} in 3--space each of the
components/strands of which is made of a finite number of straight arcs endowed with matching
orientation; also  their {\it diagrams}, that is, regular projections on the plane (with  only
finitely many crossings), where in addition no vertex of the link/braid should be mapped onto a
double point. 
\bigbreak 
\noindent A {\it geometric braid diagram} or simply a {\it braid} has a top-to-bottom 
direction and, in addition, we require that no two crossings are on the same horizontal level. 
 If we slice up in general position (i.e. without cutting through crossings) a braid diagram on
$n$ strands,  it may be seen as a word on the (well-known) basic crossings \ $\sigma_i$ \ and \
${\sigma_i}^{-1} $ \ for \ $i=1,\ldots,n-1$ (we refer the reader to figure 9 for an example). The
set of braids  on  $n$  strands modulo isotopy gives rise to the braid group $B_n$ with 
presentation: \[ B_n= \langle \sigma_1,...,\sigma_{n-1}\ |\ \sigma_i \sigma_j = \sigma_j \sigma_i 
\mbox{ \  for \ } |i-j|>1,\ \ \sigma_i \sigma_{i+1} \sigma_i=\sigma_{i+1} \sigma_i
\sigma_{i+1}\rangle .\] 
 The operation in the group is  {\it concatenation}  (we place one
braid on top of the other), and the identity element is the braid on $n$ strands that does no
braiding. Note that the elements of $B_n$ are also called `braids', but in here there will 
be no ambiguity, since we shall be working with geometric braid diagrams.
\bigbreak 
There are two combinatorial moves on diagrams which we shall consider.
\begin{itemize}
\item[(1)] {\bf $\Delta$--move:} An arc is replaced by two arcs forming a 
triangle (and its inverse), respecting orientation and crossings. 
`Respecting crossings' means that, if we lift the diagram to an embedding
in 3--space, then the $\Delta$--move lifts to an elementary isotopy (see figure 3 for examples).
\item[(2)] {\bf Subdivision:} A vertex is introduced/deleted in an arc of the diagram.
\end{itemize}

\noindent Subdivision moves may be viewed as special cases of $\Delta$--moves. We shall call the
equivalence relation generated by these two moves a {\it combinatorial isotopy} or just an {\it
isotopy}. It is a classical result of combinatorial topology that this notion of isotopy is
equivalent to the standard definition of combinatorial (or PL) isotopy of the
embedding obtained by lifting to 3--space, and this in turn is equivalent to the notion of
isotopy in the smooth category. Moreover, Reidemeister  \cite{Rd1} (and Alexander, Briggs  
\cite{AB})  proved that a $\Delta$--move can break into a finite sequence of planar 
$\Delta$--moves and the three local $\Delta$--moves illustrated below (known as `Reidemeister
moves') with their obvious symmetries (for a detailed account see \cite{BZ}).

\bigbreak
$$\vbox{\picill5.5inby1in(fig3)  }$$

\begin{defn}[$L$--moves]{\rm \ Let $D$ be a link diagram/braid  and $P$ a point of an arc of $D$
such that $P$ is not vertically aligned with any of the  crossings or (other) vertices of $D$
(note that $P$ itself may be a vertex). Then we can perform the following operation:  Cut the arc
at $P$, bend the two resulting smaller arcs apart slightly by a small isotopy and introduce
two new vertical arcs to new top and bottom end-points in the same vertical line as $P$.  The
 new  arcs are both oriented downwards and they  run either both {\it under}
or both {\it over} all other arcs of the diagram.  Thus there are two types of $L$--moves, an
{\it under $L$--move} or {\it $L_u$--move} and an {\it over $L$--move} or {\it $L_o$--move}.
(Recall figure 1 above for an example of the two braid $L$--moves.)}\end{defn}

Below we illustrate $L_o$--moves applied to an edge $QR$ of an oriented link diagram such that $QR$
moves via isotopy from a downward arc to an upward arc through the horizontal position. The thick
circle will represent `the rest of the diagram', while the region inside the circle shall be called
`the magnified region'. Note that, if we join freely (using the dotted arcs) the two new arcs we
obtain a link diagram isotopic to the one on which the $L$--move was applied.

\bigbreak
$$\vbox{\picill4.5inby1.6in(fig4)    }$$

\noindent Definition 2.1 could be given even more generally so as to accommodate the possibility
of the new vertical strands running {\it upwards} (figure 4b).
Although, for the purposes of this paper we shall only consider $L$--moves with the new vertical
arcs oriented {\it downwards}.

\begin{rem} \rm Using a small braid isotopy, a braid $L$--move can be equivalently seen with
a crossing (positive or negative) formed:

\bigbreak
$$\vbox{\picill5inby1.5in(fig5)   }$$

\noindent This gives the following algebraic expression for an $L_o$--move and an $L_u$--move
respectively.  
\[ \alpha=\alpha_1\alpha_2 \sim
\sigma_i^{-1}\ldots \sigma_n^{-1} \widetilde{\alpha_1}\sigma_{i-1}^{-1}\ldots
\sigma_{n-1}^{-1}\sigma_n^{\pm 1}\sigma_{n-1} \ldots \sigma_i
\widetilde{\alpha_2}\sigma_n \ldots \sigma_i \]  
                     \[ \alpha=\alpha_1\alpha_2 \sim
\sigma_i\ldots \sigma_n \widetilde{\alpha_1}\sigma_{i-1}\ldots
\sigma_{n-1}\sigma_n^{\pm 1}\sigma_{n-1}^{-1}\ldots\sigma_i^{-1}
\widetilde{\alpha_2}\sigma_n^{-1}\ldots\sigma_i^{-1} \] 
\noindent  where $\alpha_1$, 
 $\alpha_2$ are elements of $B_n$ and $ \widetilde{\alpha_1}$,
 $\widetilde{\alpha_2} \in B_{n+1}$ are obtained from $\alpha_1$,
 $\alpha_2$ by replacing each $\sigma_j$ by $\sigma_{j+1}$ for
 $j=i,\ldots,n-1$. \end{rem}

$L$--moves and isotopy generate an equivalence relation in the set of braids, which we shall call
{\it $L$--equivalence}. Let now $B$ be a braid.  The {\it closure} of $B$ is the oriented 
link diagram $\C(B)$ obtained by joining each top end-point to the corresponding bottom end-point
by almost vertical arcs as illustrated in figure 6, where $B$ is  contained in a `box'. (Note that
we draw some smooth arcs for convenience.)

\bigbreak

$$\vbox{\picill3.5inby1.3in(fig6)  }$$

\noindent Our aim is to prove the following theorem.

\begin{th}[One-move Markov theorem] $\C$ induces a bijection between the set of $L$--equivalence 
classes of braids and the set of isotopy types of (oriented) link diagrams.
\end{th}

\section{The braiding process}

We shall define an inverse bijection to $\C$ by means of a canonical {\it
braiding process} which turns an oriented link diagram (with a little extra
structure) into a braid. Note that we only work with oriented diagrams, so in the sequel we shall
drop the adjective. Let $D$ be a link diagram with no horizontal arcs, and consider the
arcs in $D$ which slope upwards with respect to their orientations; call these arcs {\it  opposite
arcs}. In order to obtain a braid from that diagram we want:  
\bigbreak
\noindent 1) to keep the arcs that go downwards.

\noindent 2) to eliminate the opposite arcs and produce instead  braid strands. 
\bigbreak
\noindent If we run along an opposite arc we are likely to meet a succession of
overcrossings and undercrossings. We subdivide (marking with points) every
 opposite arc into smaller -- if necessary -- pieces, each containing crossings
of only  {\it one}  type;  i.e. we may have:

\bigbreak

$$\vbox{\picill4.5inby1.2in(fig7)  }$$

\noindent We call the resulting pieces {\it  up--arcs}, and we label every 
 up--arc with an  {\it `o'/`u'} according as it is the  {\it over/under}  arc of a
crossing (or some crossings). If it is a  {\it free up--arc}
(and therefore  it contains no crossings), then we have a choice whether  to label it {\it `o'} 
or  {\it `u'}. The idea is to eliminate the  opposite arcs by eliminating their up--arcs one
by one and create braid strands instead.  Let now $P_1,P_2,\ldots,P_n$ be the top vertices of the
up--arcs; fix attention on one particular top vertex $P=P_i$ and suppose
that $P$ is the top vertex of the up--arc $QP$. 
\bigbreak
\noindent Associated to $QP$ is the {\it sliding  
triangle} $T(P)$, which is a special case of a triangle
needed for a $\Delta$--move; it is right-angled with hypotenuse  $QP$
and with the right angle lying below the up--arc. Note that,  if $QP$ is itself vertical, then
$T(P)$  degenerates into the arc $QP$.  We say that a sliding triangle is of type {\it over} or
{\it under} according to the label of the up--arc it is associated with. (This implies that there
may be triangles of the same type lying one on top of the other.)

 The germ of our braiding process is this. Suppose for definiteness that $QP$  is of
type  {\it over}. Then perform an $L_o$--move at $P$ followed by a $\Delta$--move across the
sliding triangle $T(P)$ (see figure 8). By general position the resulting diagram will be regular
and $QQ'$ may be assumed to slope slightly downwards. If $QP$ were {\it under} then the $L_o$--move
 would be replaced by an $L_u$--move. Note that the effect of these two moves has been  to
replace the up--arc $QP$ by three arcs none of which are up--arcs, and therefore
we now have fewer up--arcs. If we repeat this process for each
up--arc in turn, then the result will be a braid.

\bigbreak

$$\vbox{\picill4.5inby1.5in(fig8)  }$$ 

\noindent In the following example we illustrate the braiding process applied to a particular link
diagram (we do not draw the sliding triangles). The numbering of the strands indicates the
order in which we eliminated the up--arcs. Obviously, for different orders we
may  obtain different braids.

\bigbreak

$$\vbox{\picill5.5inby2.7in(fig9)  }$$

\noindent It is clear from the above that different choices when doing the braiding (e.g. when
choosing labels for the free up--arcs or when choosing an order for layering sliding
triangles of the same type)  as well as slight isotopy changes on the diagram level may result in
important 
 changes in the braid picture. For the proof of Theorem 2.3 we would like to control such changes as
much as possible.  For this we shall need a more rigorous setting. 
\bigbreak
Sliding triangles are said to be {\it adjacent} if the corresponding up--arcs have a common
vertex (then the sliding  triangles will have a common corner).
\bigbreak
\noindent {\bf Triangle condition} {\it Non-adjacent sliding triangles are only allowed to
meet if they are of opposite types (i.e. one over and the other under).}

\bigbreak

$$\vbox{\picill5inby1.3in(fig10)  }$$ 

 \begin{lem}{Given a link diagram $D$, there is a subdivision $D'$ of $D$ such that 
(for appropriate choices of under/over  for free up--arcs) the triangle condition is satisfied.}
\end{lem} 

\begin{pf} Let  $d =$ minimum distance between any two crossings of $D$.
Let  $0<r<d/2$  be such that any circle of radius  $r$  centred
at a crossing point does not intersect with any other arc of the diagram. Let  $s$ 
be the minimum distance between any two points in $D$ further than 
$r$  away from any crossing point. Now let  $\varepsilon= \frac{1}{2}\, $min\{$s$,$r$\}  and 
$D'$  be a subdivision of the diagram such that the length of every up--arc 
is less than  $\varepsilon$. Then the triangle condition is satisfied, provided we
make the right choices (under/over) for sliding triangles of free
up--arcs near crossings (see picture above).   \end{pf}

\begin{rems} \rm (1)  The triangle condition  implies that the
eliminating moves do not interfere with each other, so there are no pairs of sliding triangles
that need layering and therefore it does not matter  in what order we eliminate the up--arcs. In
fact we can eliminate all of them
 {\it simultaneously}. 
\bigbreak
\noindent (2)  The choice of $\varepsilon$ in the last proof is far 
smaller than necessary.  In fact we only need to consider 
  the subdiagram consisting of up--arcs and crossings
on up--arcs.
\bigbreak
\noindent (3)  Notice that the choice of under/over for some of the
free up--arcs may be forced by the labelling of other up--arcs (see for example the
top  {\it over} up--arc in figure 10).  We can -- if we wish -- further
subdivide the diagram (here the under up--arc) to make such pairs disjoint. 
\bigbreak
\noindent (4)  If a diagram satisfies the triangle condition, then so does
any subdivision (where the smaller triangles may need to retain the same labels).  \end{rems}

\begin{defn} \rm  A {\it generic diagram} is a link diagram with subdividing points and
sliding triangles put in  general position with respect to the
height function,  such that the following conditions hold: 
\bigbreak
\noindent 1) \ there are no horizontal arcs,

\noindent 2) \ no two disjoint subdividing points are in vertical alignment, where by  
`disjoint' we mean subdividing points that do not share a common edge.

\noindent 3) \ any two non-adjacent sliding triangles satisfy the triangle condition
and if they intersect, this should be along a common interior (and not a single
point). \end{defn}
  
\noindent Conditions 1 and 3 are related to the braiding, whilst condition 2 
ensures that no pair of strands in the resulting braid will be in the same vertical line. 

\begin{defn} \rm A {\it generic $\Delta$--move} is a $\Delta$--move between generic diagrams. 
\end{defn}
\noindent  In the sequel, by `$\Delta$--moves' we shall always refer to the local planar
$\Delta$--moves together with the Reidemeister moves.

\begin{lem}{An isotopy between generic link diagrams can be
realized using only generic    $\Delta$--moves. } \end{lem}

\begin{pf} If after some  $\Delta$--move during the isotopy appears a horizontal arc or vertical
alignment of vertices we remove it by replacing one of the participating vertices by a point
arbitrarily close to it, so that the new point will not cause such a singularity or violation of
condition 3 in {\it all} the diagrams of the isotopy chain. This is possible by a general position
argument. If two non-adjacent sliding
triangles touch so that condition 3 is violated (figure 11a) we argue as above. Also, by
definition of  regular isotopy, two non-adjacent triangles cannot touch on two subdividing points
or on a point with a hypotenuse or along their hypotenuses.  Therefore, the remaining possibilities
are the ones illustrated in figure 11b.

\bigbreak

$$\vbox{\picill5inby1.1in(fig11)  }$$ 

\noindent As before, by general position we can replace one of the participating vertices by a
point arbitrarily close to it, so that the new point will not violate conditions 1, 2 and 3 in 
all the diagrams of the isotopy chain. Finally, if two triangles happen to intersect after  a
$\Delta$--move and they are of opposite type, then the move is generic. If they are of the
same type and the intersection is not essential (see figure 12) we simply subdivide further and we
carry the subdivision through in the whole chain of the isotopic diagrams, taking care that the
subdividing point will not violate condition 2 in the whole chain. (Note that, by Remark 3.2
(4), condition 3 will not be violated.)  

\bigbreak

$$\vbox{\picill5inby1.4in(fig12)  }$$

\noindent  If the intersection is essential (for example in figure 13 a free up--arc labelled `u'
moves by a $\Delta$--move over another arc), then we introduce an appropriate extra subdividing
point so as to create a smaller free up--arc to which we attach the opposite label. This is always
possible by Lemma 3.1.

\bigbreak

$$\vbox{\picill2.8inby1.1in(fig13)  }$$

\end{pf}

\noindent It is clear from Lemma 3.5 that generic link diagrams are dense in the space of all
diagrams;  indeed, a non-generic link diagram may be seen as a middle stage of the isotopy between
generic diagrams.  Thus w.l.o.g. we shall assume from now on that all diagrams are generic; also,
by virtue of Lemma 3.5 that all  isotopy moves are also generic.

We are now ready to give a rigorous braiding process. Namely, take a link diagram and eliminate
one by one the up--arcs in the way described above. By Remark 3.2 (1) the order of the eliminating
moves is now irrelevant.

\begin{cor}[Alexander's theorem]  Any (oriented) link diagram is  isotopic to the closure of a braid
\footnote{H. Brunn \cite{Br} in 1897 proved that any link has a projection with a single
multiple point; from which it follows immediately (by appropriate perturbations) that we can braid
any link diagram. Other proofs of Alexander's theorem have been given by  H.R. Morton \cite{M},
 S. Yamada \cite{Y}, P. Vogel \cite{V}.}. \end{cor}

\begin{pf} The braiding process comprises $L$--moves and then isotopies.  But 
the effect of an $L$--move after closure is by definition an isotopy. \end{pf}

\section{Proof of Theorem 2.3}

By the local nature of the $\Delta$--moves we may assume that for a given  link diagram 
we have done the braiding for all up--arcs except for the ones that we are interested in every
time; these will be lying in the magnified region placed inside the braid. By Remark 3.2 (1) this
choice does not affect the final braid.
\bigbreak
Now, two braids that differ by a finite sequence of $L$--moves  have isotopic closures.
Therefore the function $\C$ from  $L$--equivalence classes of braids to isotopy types of link
diagrams is well-defined. To show that $\C$ is a bijection we shall use  our braiding process to
define an inverse function $\B$. Namely, for a diagram $D$ let $\B(D)$ be the braid resulting from
the braiding algorithm applied to it.  We have to show that $\B$ is a 
well-defined function from link-diagram types to $L$--equivalence classes of braids, therefore we
have to check that $\B(D)$ does not depend up to $L$--equivalence on  the choices made before the
braiding and on  $\Delta$--moves between link diagrams.  
\bigbreak
\noindent The choices made before the braiding
consist of the subdividing points we choose and the labelling we may have to choose for some free 
up--arcs.
 Finally, we have to show that  $\C$ and $\B$ are
mutually inverse. This is easy: Closing the result of the braiding process yields a link isotopic
to the original one, therefore $\C\circ \B = \text{id}$. Moreover, applying the braiding process to
the closure of a braid yields precisely the same braid back again (from the way we defined closure),
so $\B\circ \C = \text{id}$.

The proof relies entirely on Remark 3.2 (1)  and on the following two lemmas.  
  
\begin{lem}{If we add on an up--arc, $\alpha$, an extra subdividing point  $P$  and label the two
new up--arcs, $\alpha_1$  and  $\alpha_2$, the same as $\alpha$, the
corresponding braids are $L$--equivalent.} \end{lem}

\begin{pf} For definiteness we assume that $\alpha$ is labelled with an `o'. We
complete the braiding of the original diagram by eliminating  $\alpha$  (see picture below). Then,
on the new horizontal piece of string, we take  an arbitrarily small neighbourhood  $N'$  around 
$P'$, the projection of $P$ (see picture below). By general position $N'$  slopes slightly
downwards. We then perform an {\it over}  L-move at $P'$. Finally, sliding an appropriate piece  
of string using braid planar isotopy, we obtain the braid that would result from the original
diagram with the subdividing point  $P$  included (see figure 14). (Note  that new vertical
strands coming from the braiding that may run over or under $\alpha$  do not affect the result.)

\bigbreak

$$\vbox{\picill5.5inby1.9in(fig14)  }$$

\end{pf}
\begin{lem}{When we meet a free up--arc, which we have the choice of
labelling  `u'  or `o', the resulting braid does not depend  -- up
to $L$--equivalence -- on this choice.} \end{lem}
 
\begin{pf} First, we shall assume for simplicity that the
sliding triangle of the up--arc does not lie over or under any other arcs of the
original diagram. Also we assume for definiteness that the up--arc is
originally labelled  `o'. We complete the braiding by eliminating it. Then on
the new  almost-horizontal piece of string we take an arbitrarily small neighbourhood
 $N'$ of a point $P'$ that is a projection of an arbitrarily small neighbourhood $N$ on the
original up--arc, and such that  there is no other vertical line between the vertical line of $P$
and the one of $P'$ (figure 15).  We then perform an $L_u$--move  at  $P'$. 

\bigbreak

$$\vbox{\picill5inby1.9in(fig15) }$$

\noindent The fact that the original up--arc is free and small enough implies
that only vertical strands can pass over or under its sliding
triangle. Therefore -- as  $N$  is arbitrarily small -- there is no arc
crossing  $AB$  so as to force it be an  {\it under}  arc. Also, by braid
planar isotopy we shift  $A$  slightly  {\it higher},
 so as to come to the position where we can undo an $L_o$--move (see figure 16). We undo it, so the
final braid -- up to a small braid planar isotopy -- can be seen as the braid that we would have
obtained from the original diagram  with the free up--arc labelled with  `u'  instead of  `o'.

\bigbreak

$$\vbox{\picill6inby1.9in(fig16) }$$

\noindent Notice that, if the original up--arc were an  `u' we would perform an 
$L_o$--move  at $P'$. To complete the proof of the lemma we assume that the sliding
triangle of our up--arc lies over or under other arcs of the original diagram.
In this case we subdivide it (using Lemma 4.1) into arcs small enough to ensure
that all the sliding triangles are clear; we give all new arcs the labelling
of the original one. Then we change the labelling of each up--arc using the
above and, using Lemma 4.1 again, we eliminate all the new subdividing points (figure 17). 

\bigbreak

$$\vbox{\picill5inby1.4in(fig17) }$$

\end{pf}

\begin{cor}{ If we have a chain of overlapping sliding triangles of free
up--arcs so that we have a free choice of labelling for the whole chain  
 then, by Lemmas 4.1 and 4.2, this choice does not affect  -- up to
$L$--equivalence -- the final braid.} \end{cor}

\begin{cor}{If by adding a subdividing point on an up--arc we have a choice
for relabelling the resulting new up--arcs  so that the triangle condition  is
still satisfied then, by Lemmas 4.1 and 4.2, the resulting braids are
$L$--equivalent.} \end{cor}

\begin{cor}{Given any two subdivisions, $S_1$ and $S_2$, of a diagram which  will satisfy the
triangle condition with appropriate  labellings, the resulting braids  are $L$--equivalent.}
\end{cor}
 
\begin{pf} By Remark 3.2 (4) this can be easily seen if we consider the subdivision
  $S_1 \bigcup S_2$  and apply the lemmas above.\end{pf}

\noindent Corollary 4.5 proves independence of subdivision and  labelling for  diagrams, and thus
we are done with the static part of the proof. 
 Using Corollary 4.5 we have also proved independence of  $\Delta$--moves related to 
 condition 3 of Definition 3.3 (recall proof of Lemma 3.5). 
\bigbreak
It remains to consider the effect of
the rest of the  $\Delta$--moves. We shall check first the planar $\Delta$--moves. These will
include the examination of $\Delta$--moves related to conditions 1 and 2 of Definition 3.3. Now, a
$\Delta$--move can be regarded as a continuous family of diagrams (see figure 18). 
From the above we shall assume that the triangle condition is not violated. Also, by symmetry, we
only have to check the $\Delta$--moves that take place in the first quarter of the plane, above the
edge $AB$, and w.l.o.g. $AB$ is an up--arc. Therefore, the $\Delta$--moves we want to check can 
 split into moves that take place above $AB$ within the vertical zone defined by $A$ and $B$ and
moves that take place outside this zone. All resulting braids will be compared with the braid
obtained by subdividing $AB$ at an appropriate point $P$. I.e. we have:

{\it (i) Inside subdivision} \  Consider the continuous family of triangles with edge $AB$  
and all the new vertices lying on the vertical line of an interior point $P$ of $AB$
(figure 18). A new vertex cannot aligne horizontally with 
$B$. So, if $P'B$ is a new edge before the horizontal position from below, we can see
$L$--equivalence of the corresponding braids by introducing an extra vertex in $AB$ at $P$  
 and applying Lemma 4.1 (see left-hand side of figure 19). 

\bigbreak

$$\vbox{\picill4.5inby1.6in(fig18)  }$$

\noindent If $P''B$ has passed above the horizontal position and therefore it becomes a
down--arc  (cf. Definition 3.3, condition 1), we introduce an  $L$--move at the point $B$
 and we show that again we obtain the same braid as if we had originally subdivided $AB$ at $P$
(see right-hand side of figure 19). The type of the  $L$--move  is over/under according as $AB$ is
an over/under up--arc.  

\bigbreak

$$\vbox{\picill5.5inby1.6in(fig19)  }$$

\noindent  Note that, if $AB$ were originally a down--arc, then the
new edges would pass from down--arcs to opposite arcs through the horizontal position. In 
this case,  according to the way we do the braiding and to the triangle condition, we may need to
subdivide further the new opposite arcs and apply lemmas 4.1 and 4.2.   
\bigbreak
\noindent  Before continuing with the proof we introduce some extra notation. We shall denote  by
$(APB)$ the braid obtained by completing the braiding on the arcs $AP$ and $PB$. Suppose now we
subdivide $AB$ at a different point, $Q$ say. Then, by Lemma 4.1 we have $(APB) \sim (AQB)$, since
they both differ by an $L$--move from $(APQB)$. In figure 20 the dotted line indicates a possible
braid strand between the vertical lines of $P$ and $Q$ (cf. Definition 3.3, condition 2). 

\bigbreak

$$\vbox{\picill2.1inby1.4in(fig20)  }$$

{\it (ii) Outside `subdivision'} \  Consider now the continuous family of triangles  with edge
$AB$  and all the new vertices lying on the vertical line of a point $P$ that lies outside the
vertical zone defined by $A$ and $B$ (see figure 21). By symmetry we only need to examine the
$\Delta$--moves taking place on the side of $B$, say.  Assume first that $P$ is close enough to $B$
so that no braid strand passes in between  the vertical lines of $P$ and $B$. Then obviously 
$(APB) \sim (AB)$ (braid isotopy). Reasoning as above, if $P'$ is a new vertex such
that $P'B$ is below the horizontal position, then $(AP'B)$ differs from $(APB)$ by an $L$--move
introduced at $B$. Also, if $P''$ is a new vertex such that $P''A$ is below the horizontal, then 
$(AP''B)$ differs from $(AP'B)$ by an $L$--move introduced on $AP''$ at the point $P''$.     

\bigbreak

$$\vbox{\picill4inby1.6in(fig21)  }$$

\noindent  Assume now that the $\Delta$--move
introduces a new vertex $Q$ so that between the vertical lines of $P$ and $Q$  there is a braid
strand (cf. Definition 3.3, condition 2). In this case we see the $L$--equivalence by introducing
and deleting two $L_o$--moves (figure 22).   

\bigbreak

$$\vbox{\picill5.5inby2in(fig22)  }$$

We shall now check the Reidemeister moves (recall figure 3).  
\bigbreak
\noindent To check the first two moves we follow similar
reasoning as above, but we shall demonstrate it here for completeness. For the first one, we 
illustrate below that the braid obtained after the performance of the  move is equivalent to the
one obtained before the move, up to  braid isotopy and one
$L_o$--move:

\bigbreak

$$\vbox{\picill6inby1.8in(fig23)  }$$

\noindent For the second move we complete the braiding of the left-hand side diagram and we
notice that, using braid isotopy, we can undo  an $L_u$--move. We then
obtain the braid that we would obtain from the right-hand side diagram with label  `u'  for the
up--arc.

\bigbreak

$$\vbox{\picill4inby1.8in(fig24)  }$$    

We shall now check the third type of Reidemeister moves or `triple point moves'. If all three
arcs are down--arcs then the move is a braid isotopy. If either of the outer arcs is an up--arc
then the move is invisible in the braid. It remains to check the case where both outer arcs are
down--arcs and the middle arc is an up--arc. But then the following trick changes the situation
to the one where all three arcs are down--arcs. 

\bigbreak

$$\vbox{\picill5inby1.1in(fig25)  }$$

We have now proved that the function $\B$ from isotopy types of link diagrams to
$L$--equivalence classes of braids is well-defined, and the proof of Theorem 2.3 is now completed.  
\ \hfill $\Box$

\subsection{A comment on conjugation} 

The classical Markov braid theorem states that: {\it Isotopy classes of (oriented) link diagrams 
are in 1-1 correspondence  with  certain equivalence classes in the set of all braids, the  
equivalence  being given by the following two algebraically formulated moves between braids in 
$\bigcup_{n=1}^{\infty}B_n$: 
\begin{itemize}
\item[(i)]{\it Conjugation:} If \ $\alpha,\beta \in B_n$
 \ then \ $\alpha \sim \beta^{-1} \alpha \beta$. 

\item[(ii)]{\it Markov moves or $M$--moves:} If \ 
 $\alpha \in B_n$ \ then \ $\alpha \sim \alpha\sigma_n^{+1}\in B_{n+1}$ \ and 
 \ $\alpha\sim \alpha\sigma_n^{-1} \in B_{n+1}$. 
\end{itemize} }

 It is clear from Remark 2.2 that an $M$--move is a special case of an $L$--move.
Also, it follows from the proof of Theorem 2.3 that conjugation can be realized by a (finite)
sequence of $L$--moves. Indeed, let $\alpha$ and $\beta ^{-1}\alpha\beta$ be two conjugate braids.
Then their closures  $\C(\alpha)$ and $\C(\beta ^{-1}\alpha\beta)$ are isotopic diagrams and so
they differ by a sequence of $n$, say, $\Delta$--moves; so, for each stage of
isotopy we obtain the following sequence of link diagrams 
\[ \C(\beta ^{-1}\alpha\beta) \sim L_1 \cdots \sim L_{n-1} \sim \C(\alpha) \]
which we turn into the braids $B_0, B_1, \cdots, B_{n-1}, B_n$, say, using our algorithm. From
the above, every consecutive pair of these braids differs by a finite sequence of $L$--moves. Now, 
 the way we defined the closure of a braid (recall figure 6) guarantees that $B_0$
is  $\beta ^{-1}\alpha\beta$ and $B_n$ is $\alpha$. Thus, conjugation is only a redundant
`auxilliary' move used to bring the $M$--moves inside the braid box, and therefore Theorem 2.3
indeed sharpens the classical result.
 
As a concrete example, we illustrate below the main instances of the isotopy sequence from
the closure of a braid that is conjugate to $\alpha$ by an elementary crossing, up to
$\C(\alpha)$.   

\bigbreak

$$\vbox{\picill4inby1.6in(fig26)  }$$ 

\begin{rems} \rm (1) The proof of Theorem 2.3 may be clearly used as an alternative proof of the
classical theorem. Indeed, whenever in the proof appears an $L$--move we use conjugation to bring
it to the position of an $M$--move. The same reasoning holds for all Markov-type theorems that
follow in the sequel. 
\bigbreak
\noindent (2) It should be stressed that Theorem 2.3 sharpens the classical result only if 
we work with  open braids. If we work with closed braids we do not need conjugation in the braid
equivalence. We should also stress that vertical alignment of top vertices is the most
important  case to check in the proof of Theorem 2.3 because it corresponds precisely to
conjugation by a generator $\sigma_i$ in $B_n$. \end{rems}   
 
\subsection{Relative version of Markov's theorem}

As a consequence of the proof of Theorem 2.3 we can prove a relative version of the result. By a
{\it braided portion} of a link diagram we mean a finite number of arcs in the link which are all
oriented downwards (so that the braid resulting from our braiding process will contain this
braided portion).

\begin{th} Let $L_1$, $L_2$ be oriented link diagrams which both contain a common braided portion
$B$. Suppose that there is an isotopy of $L_1$ to $L_2$ which finishes with a homeomorphism fixed
on $B$. Suppose further that $B_1$ and $B_2$ are braids obtained from our braiding process applied
to $L_1$ and $L_2$ respectively. Then $B_1$ and $B_2$ are $L$--equivalent by moves that do not
affect the common braided portion $B$. \end{th} 

\begin{pf} Assume first that the isotopy from $L_1$ to $L_2$ keeps $B$ fixed; the result then
follows immediately from the proof of Theorem 2.3 because the braided portion does not participate
in the proof. We shall use some standard PL topology to reduce the general case to this special 
case. 
\bigbreak
Let $N$ denote a relative regular neighbourhood of $B$ relative to the ends of the arcs
(i.e. $N$ comprises a number $t$ say of 3--balls each of which contains an arc of $B$ as an
unknotted subarc). Think momentarily of these balls as small balls centred at points $P_1,
\ldots, P_t$. The isotopy restricted to $P_1, \ldots, P_t$ determines a loop in the configuration
space of $t$ points in ${\bf R}^3$. But this configuration space is well-known to be
simply-connected and therefore we may assume that the isotopy is fixed on $P_1, \ldots, P_t$. By
the  regular neighbourhood theorem we may now assume that the isotopy fixes $N$ setwise. 
\bigbreak
\noindent Now
restrict attention to one of the balls. The isotopy on this ball determines a loop in the space of
the  PL homeomorphisms of the 3--ball. But $\pi_1$ of this space is generated by a rotation through
$2\pi$ about some axis. We may suppose that this axis is the corresponding unknotted arc of $B$ and
hence that the isotopy fixes this subarc pointwise. Thus we may assume that the whole of $B$ is
fixed pointwise. The result now follows from the special case. 
\end{pf}

\section{Extension of results to other 3--manifolds}
\def\sectionname{Extension to 3--manifolds}

In this section we  use the methods of Theorem 2.3 to prove a second relative version (relative to
a fixed closed subbraid) and deduce an analogue of Markov's theorem for isotopy of oriented links in
knot/link complements in $S^3$.   This is then used for extending the results to closed connected
orientable 3--manifolds. 

\subsection{Markov's theorem in knot complements}

Let $S^3 \backslash K$ be the  complement of the oriented knot $K$ in $S^3$. By `knot
complement' we  refer to complements of both knots and links. Using  Alexander's theorem and the
definition of ambient isotopy, $S^3 \backslash K$ is homeomorphic to $M=S^3 \backslash
\widehat{B}$,  where $\widehat{B}$  is isotopic to $K$ and it is the closure of some braid $B$. Let
now $L$ be an oriented link in $M$. If we fix $\widehat{B}$ pointwise on its projection plane we
may represent $L$ unambiguously  by the {\it mixed link} $\widehat{B}\bigcup L$ in $S^3$, that is,
a link in $S^3$  consisting of the {\it fixed part} $\widehat{B}$ and the {\it standard part} $L$
that links with $\widehat{B}$  (for an example see figure 27a). 

\begin{defn}{\rm A {\it mixed link diagram} is a diagram $\widehat{B} \bigcup
\widetilde{L}$ of $\widehat{B}\bigcup L$ projected on the plane of $\widehat{B}$ which is 
equipped with the top-to-bottom direction of $B$.}\end{defn} 

\bigbreak

$$\vbox{\picill5.5inby1.5in(fig27)  }$$ 

\noindent Let now $L_1$, $L_2$ be two oriented links in $M$. It follows from standard results  of
PL Topology that  $L_1$ and $L_2$ are isotopic in $M$  if and only if the mixed links
$\widehat{B}\bigcup L_1$ and $\widehat{B}\bigcup L_2$ are isotopic in $S^3$ by an ambient
isotopy which keeps $\widehat{B}$ pointwise fixed. See for example \cite{RS}; chapter 4. In terms
of mixed diagrams this isotopy will involve the fixed part of the mixed links only in the moves
illustrated in figure 27b. These shall be called `extended Reidemeister moves'. Therefore 
Reidemeister's theorem generalizes as follows in terms of mixed link diagrams.

\begin{th} Two (oriented) links in $M$ are isotopic if and only if any two corresponding mixed link
diagrams (in  $S^3$) differ by a finite sequence of the extented Reidemeister
moves together with the planar $\Delta$--moves and the Reidemeister moves for the standard
parts of the mixed links. \end{th} 

We now wish to find an analogue of Markov's theorem for links in $M$. Since  $\widehat{B}$ must
remain fixed  we need a braiding process for mixed link diagrams that  maps  $\widehat{B}$ to $B$
(and {\it not} an $L$--equivalent braid or a conjugate of $B$). Such braids shall be called {\it
mixed braids} (see figure 28c for an example).  Also the subbraid of a mixed braid that
complements  $B$ shall be called {\it the permutation braid}.
 
\begin{th}[Alexander's theorem for $S^3 \backslash \widehat B$] \ 
\ Any (oriented) link in $S^3 \backslash \widehat B$  can be represented  in $S^3$ by
some mixed braid  $B_1\bigcup B$, the closure of which is isotopic to a mixed link diagram 
$\widetilde{L_1} \bigcup\widehat B$ representing the link.
 \end{th}
\begin{pf} The fixed part  $\widehat{B}$ may be viewed in $S^3$ as the braid $B$ union an arc, $k$
say, at infinity. This arc is the identification of the two horizontal arcs containing the endpoints
of $B$ and thus it realizes the closure of $B$. (In figure 28a the braid $B$ is drawn curved but
this is just for the purpose of picturing it.) Let $N(k)$ be a small regular neighbourhood of $k$
and let $L$ be a link in $M$. By general position $L$ misses  $N(k)$ and
therefore it can be isotoped into the complement of $N(k)$ in $S^3$. 

By expanding $N(k)$ we can
view its complement as the cylinder $T=D^2\times I$ that contains $B$ (figure 28b). We then
apply our braiding in $T$. This will leave $B$ untouched but it will braid $L$ in $T$. 
Finally we cut $\widehat{B}$ open along $k$ in order to obtain a mixed braid (figure 28c).
\end{pf}

\bigbreak

$$\vbox{\picill5.5inby2.2in(fig28)  }$$

\begin{note} \rm An alternative more `rigid' proof modifying the braiding of section 3 is given in
\cite{La},\cite{La1}\footnote{Theorem 5.3 can be also proved using for example  Alexander's original
braiding  (\cite{AB}, \cite{B}) or alternatively Morton's threading (\cite{M}), both resulting
closed braids, so that the braid axis/thread goes through a point on the plane around which 
all strings of $\widehat{B}$ have counter-clockwise orientation; then we cut the braid open along
a line that starts from the central point and  cuts through the closing side of
$\widehat{B}$.}. \end{note}

\noindent We are now ready  to prove the second relative version of Markov's theorem and a
version for knot complements. 

\begin{th} \ Two  oriented links in $S^3$ which contain a common subbraid $B$ are isotopic rel
$B$  if and only if any two corresponding mixed braids $B_1\bigcup B$ and $B_2\bigcup
B$  are $L$--equivalent by moves that do not touch the common subbraid $B$. 

\noindent Further  any
oriented links  in $S^3 \backslash \widehat B$ are isotopic  in $S^3 \backslash \widehat B$  if and
only if any two corresponding mixed braids $B_1\bigcup B$ and $B_2\bigcup B$ in $S^3$ are
 $L$--equivalent by moves that do not touch the common subbraid $B$.   \end{th} 

\begin{pf} We prove the relative version. The version for knot complements then follows
immediately from the discussion above. 

Consider an isotopy of a link $L$ in the complement of
$\widehat B$. By a general position argument we may assume that this isotopy only crosses the arc
$k$ a finite number of times and that these crossings are clear vertical cuts (figure 29a). The
neighbourhood $N(k)$ misses also the isotopy apart from those clear crossings. As in the proof
above, we expand $N(k)$ and thus the isotopy now takes place inside $T$ apart from the crossings
with $k$. Applying now the proof of Theorem 2.3, the isotopies within $T$
can be converted into a sequence of $L$--moves that do not affect $B$. Therefore it only  remains
to  check $L$--equivalence at  a crossing of $L$ with $k$. But after expanding $N(k)$ such a
crossing looks like exchanging a braid strand, $l$ say, that runs along the front of the mixed
braid for one with the same endpoints that runs along the back of the braid (figure 29b). 

\bigbreak

$$\vbox{\picill5.5inby1.8in(fig29)   }$$

\noindent If $l$ is an up--arc it can be made into a free up--arc and so it dissappears after the
braiding. If  $l$ is a down--arc we introduce an $L_u$--move at its top part so that
 between the endpoints of the new strand and the endpoints of $l$ there is no other braid
strand. We then undo an  $L_o$--move as depicted below, and the proof is concluded. Note
that it is necessary to use the version of an  $L$--move with a crossing (recall figure 5).
  \end{pf}

\bigbreak

$$\vbox{\picill5.3inby1.9in(fig30)   }$$

\subsection{Markov's theorem in oriented 3--manifolds}

Let $M$ be a closed connected orientable (c.c.o.) 3--manifold.  It is well--known \cite{L},
\cite{Wa}, \cite{R} that $M$ can be obtained by surgery on a framed link in $S^3$ with integral
framings and w.l.o.g. this link is the closure $\widehat B$ of a braid $B$. We shall refer to
$\widehat B$ as the {\it surgery link} and we shall write $M=\chi (S^3, \widehat B)$. Moreover by 
 the proof given in  \cite{L} it can be  assumed that all the components of the surgery link are
unknotted and furthermore that they form the closure of a pure braid. Thus we may assume that
$B$ is a  pure braid. So $M$ may be  represented in $S^3$ by the framed $\widehat B$ and if we fix
$\widehat B$ pointwise we have that links/braids in $M$ can be unambiguously represented by mixed
links/braids in $S^3$ -- exactly as discussed in 5.1. Thus Theorem 5.3 holds also in this case.
Only here, the braid $B$ is framed and we refer to it as the {\it surgery braid}. 
\bigbreak
Now a link $L$ in $M$ may be seen as a link in $S^3 \backslash \widehat B$ with the extra
freedom to slide across the 2--discs bounded in $M$ by the specified longitudes of the components 
of $\widehat B$. 

\begin{defn} \rm  Let  $b$ be the oriented boundary of a ribbon  and let $L_1 \bigcup \widehat B$
and $L_2 \bigcup \widehat B$ be two oriented mixed links, so that $L_2 \bigcup \widehat B$ is
the band connected sum (over $b$) of a component, $c$, of $L_1$ and the specified (from the
framing) longitude of a surgery component of $\widehat B$.  This is a non--isotopy move in $S^3$
that reflects isotopy between  $L_1$ and  $L_2$ in $M$ and we shall call it {\it  band
move}. 
\bigbreak
\noindent  A  band move can be thus split in two steps: firstly, one of the small edges of $b$ is
glued to a part of $c$ so that the orientation of the band agrees with the orientation of  $c$. The
other small edge of $b$, which we shall call {\it little band} (in ambiguity with the notion of a
band), approaches a surgery component of $\widehat B$ in an arbitrary way. Secondly, the little
band is replaced by a string running in parallel with the specified longitude of the surgery
component in such a way that the orientation of the string agrees with the orientation of $b$ and
the resulting link is $L_2 \bigcup \widehat B$.  \end{defn}

 Since we consider oriented links in $M$ there are two types, $\alpha$ and $\beta$ say, of
band moves according as in the second step the orientation of the string replacing
the little band  agrees (type  $\alpha$)  or disagrees (type  $\beta$)  with the orientation of
the surgery component -- and implicitely of its specified longitude. (See figure 31, where $p$ is an
integral framing of the surgery component.)

\bigbreak

$$\vbox{\picill5.4inby1.2in(fig31a)   }$$

$$\vbox{\picill5.4inby1.3in(fig31b)   }$$

\noindent The two types are related in the following sense.

\bigbreak

$$\vbox{\picill5.7inby1.3in(fig32)   }$$

\begin{rems} \rm (1)  Let $L_1 \bigcup \widehat B$ and $L_2 \bigcup \widehat B$ represent the  
isotopic links $L_1$ and $L_2$ in $M$. If in the isotopy sequence there is no band move involved on
the mixed link level, then $L_1$ and $L_2$ are also isotopic in $S^3 \backslash \widehat B$. In
particular, the first step of a band move reflects isotopy in $S^3 \backslash \widehat B$. So,  
 from now on whenever we say  `band move' we will always be referring to the realization of the
second step of a band move.
\bigbreak
\noindent (2) A band move, that is, the second step of the move described in Definition 5.6, takes
place in an arbitrarily thin tubular neighbourhood of the component of the surgery link that
contains no other part of the mixed link; and since $\widehat B$ is pointwise fixed it follows that
this tubular neighbourhood projects the same for any diagram of $L_2 \bigcup \widehat B$. So by
`band move' we may unambiguously refer to both the move in the 3--space and its projection. 
\end{rems}

The discussion can be summarized by saying that link isotopy in $M$ can be regarded as isotopy in 
$S^3\backslash \widehat B$ together with a finite (by transversality) number of band moves.
Therefore Theorem 5.2 extends to: 

\begin{th}[Reidemeister's theorem for $\chi (S^3, \widehat B)$] \ Two (oriented) links  in $M$
are isotopic if and only if any two corresponding mixed link diagrams differ by  a
finite sequence of the band moves, the extented Reidemeister moves  and also by  planar
$\Delta$--moves and the  Reidemeister moves for the standard parts of the mixed links.
\end{th} 

Note now that neither of the two types of band moves can appear as a move between braids; so in
order to state our extension of Markov's theorem we modify the band move of type  $\alpha$ 
appropriately by twisting the little band before performing the move. So we have the following.   

\begin{defn} \rm A {\it braid band move} or, abbreviated, a {\it  b.b.--move}  is a move between
mixed braids that reflects isotopy in $M=\chi (S^3, \widehat B)$ and is described by the following
picture (where the middle stage is only indicative).
\bigbreak

$$\vbox{\picill5.3inby1.65in(fig33)    }$$

\end{defn}

\noindent Note that a braid band move can be positive or negative depending on the type of crossing
we choose for performing it. 

\begin{th}[Markov's theorem for $M=\chi (S^3, \widehat B)$] \ Let  $L_1$, $L_2$ be two oriented
links in $M$ and let  $B_1 \bigcup B$,  $B_2 \bigcup B$  be two corresponding mixed braids in 
$S^3$. Then  $L_1$  is isotopic to $L_2$  if and only if  $B_1 \bigcup B$  is equivalent to 
$B_2 \bigcup B$ by the braid band moves and  under $L$--equivalence that does not affect $B$.
\end{th}
 
\begin{pf} Let $\widetilde{L_1}\bigcup \widehat{B}$ and $\widetilde{L_2} \bigcup  \widehat{B}$ be
two mixed link diagrams of the mixed links representing $L_1$ and $L_2$. By Theorems 5.8 
and 5.5 we only have to check that if  $\widetilde{L_1} \bigcup \widehat{B}$ and $\widetilde{L_2}
\bigcup  \widehat{B}$ differ by a band move then  $B_1 \bigcup B$ and $B_2 \bigcup B$ differ by
b.b.--moves and $L$--equivalence that does not affect $B$. In  $\widetilde{L_1} \bigcup
\widehat{B}$  the little band would be like
\bigbreak

$$\vbox{\picill5.5inby1.4in(fig34)    }$$

\noindent depending on the orientation. If the little band is an opposite arc, w.l.o.g. we may
assume that it satisfies the triangle condition. The algorithm we use ensures that we may assume
that  $\widetilde{L_1} \bigcup \widehat{B}$  and  $\widetilde{L_2} \bigcup \widehat{B}$  are braided
everywhere except for the little band in
  $\widetilde{L_1} \bigcup \widehat{B}$  (if it is an opposite arc)
and its replacement after the performance of the band move. This happens because
the band move takes place arbitrarily close to the surgery strand; so we
can produce such a zone locally in the braid (figure 35), and consequently a b.b.--move
cannot create problems with the triangle condition.

\bigbreak

$$\vbox{\picill5.5inby1.9in(fig35)   }$$

\noindent Moreover the new strand from the band move (as far as other crossings
are concerned) behaves in the same way as the surgery string itself. So,
whenever we meet other opposite arcs we label them in the same way that we
would do if the new strand were missing.

\bigbreak

$$\vbox{\picill3.4inby1.5in(fig36)  }$$ 

\noindent Hence the different cases of applying a band move to  $\widetilde{L_1}
\bigcup \widehat{B}$  amount to the following (with proofs).
\bigbreak

$$\vbox{\picill8inby1.7in(fig37)  }$$

\noindent {\bf Proof.}  We start with the left part of move  {\it (i)} and we
twist the little band (isotopy in $S^3 \backslash \widehat B$) using a {\it negative} crossing. 
Then we perform a b.b.--move with a {\it positive} crossing and we end up with the right part of
move  {\it (i)} (see figure 38).
\bigbreak

$$\vbox{\picill8inby1.6in(fig38)  }$$ 

\bigbreak

$$\vbox{\picill8inby1.6in(fig39)  }$$

 \noindent {\bf Proof.}  We start with the right part of move  {\it (ii)}. In the front of the
otherwise braided link we do a twist of the new string using a {\it negative }  crossing (see
figure 40). Then we consider the little twisted arc as a little band and we perform another
b.b.--move around the same surgery component. This second band move takes place closer to the
surgery component than the first one. Now, the shaded region in the picture below is formed by two
similar sets of opposite twists of the same string around the surgery string. So it bounds a disc
(together with the little band that is missing), the circumference of which is  {\it not} linked
with the surgery component; but this is isotopic in $S^3 \backslash \widehat B$  to the left part
of move  {\it (ii)}. I.e. move {\it (ii)} is a finite sequence of $L$--moves and b.b.--moves.

\bigbreak

$$\vbox{\picill4inby1.7in(fig40a)  }$$ 

\bigbreak

 $$\vbox{\picill5.7inby2.2in(fig40b)  }$$ 
  
\noindent Note that in the pictures above we have included another string of
the mixed braid that links with the surgery component. Clearly this does not
affect the proof. The proof of Theorem 5.10 is now concluded. \end{pf}

\begin{rems} \rm (1)   The proof holds even if the surgery braid is not a pure braid. In this case
a b.b.--move is modified so that the replacement of the little band links only with one of the
strings of the same surgery component and runs in parallel to all remaining strings of the surgery
component.

\bigbreak

 $$\vbox{\picill4inby1.6in(fig41)  }$$

\noindent (2) An analogue of Reidemeister's theorem for c.c.o. 3--manifolds is given in \cite{Su}
and an analogue of Markov's theorem for c.c.o. 3--manifolds is given in 
\cite{S}, both using the intrinsic structure of the manifold. 
\bigbreak
\noindent (3) Theorems 5.5 and  5.7 above are described only geometrically. As shown in
\cite{La},  \cite{La2} if our manifold  is the complement of the unlink or a connected sum of lens
spaces then in the set of mixed braids we have well-defined group structures, while in the other
3--manifolds we have coset structures, which enable us to give the above theorems  an algebraic
description. In particular, if our space is a solid torus (i.e. the complement of the unknot) or a
lens space $L(p,1)$ the related braid groups are the  Artin groups of type $\cal B$ (for  
details and Jones-type knot invariants we refer the reader to \cite{La}, \cite{La1}, \cite{GL}).  

\end{rems}


\bigbreak
\noindent {\sc S.L.: 16 Mill Lane, DPMMS, University of Cambridge, Cambridge CB2 1SB, U.K., E-mail
sofia@pmms.cam.ac.uk and SFB 170, `Geometrie und Analysis', Bunsenstrasse 3-5, 37073
G\"{o}ttingen, Germany.}  
\bigbreak
\noindent {\sc C.P.R.: Mathematics Institute, University of Warwick, Coventry CV4 7AL, U.K.,
E-mail cpr@maths.warwick.ac.uk}

\end{document}